\documentclass[11pt]{article}
\usepackage{amsfonts}
\usepackage{mathtools}
\usepackage{amsmath}
\usepackage{amsthm}
\usepackage{amssymb}
\usepackage{mathrsfs}
\usepackage{graphicx}
\usepackage{upgreek}
\def\pl{\partial}

\def\*{\raisebox{.5mm}{*}}

\def\div{{\,\rm div\,}}

\def\pa{\partial}

\def\R{\mathbb{R}}
\def\A{{\cal A}}
\def\H{{\cal H}}

\def\Ga{\Gamma}

\def\<{\langle}
\def\>{\rangle}

\def\Om{\Omega}
\def\OM{\Omega}
\def\Re{{\mbox{\rm Re\,}}}
\def\Im{{\mbox{\rm Im\,}}}

\def\det{\mbox{det\,}}

\def\bma{\left[\begin{array}}
\def\ema{\end{array}\right]}
\def\bda{\left|\begin{array}}
\def\eda{\end{array}\right|}
\def\be{\begin{equation}}
\def\ee{\end{equation}}

\newtheorem{thm}{{}\hskip\parindent Theorem}[section]
\newtheorem{lem}{{}\hskip\parindent Lemma}[section]
\newtheorem{pro}{{}\hskip\parindent Proposition}[section]
\newtheorem{exl}{{}\hskip\parindent Example}[section]

\newtheorem{dfn}{{}\hskip\parindent Definition}[section]
\newtheorem{rem}{{}\hskip\parindent Remark}[section]

\large
\title{Asymptotic behavior of  the nonlinear Schr\"{o}dinger  equation on  exterior domain}
\date{}
\author{Zhen-Hu Ning\thanks{Corresponding author, E-mail address: nzh41034@163.com}}
\begin{document}
\maketitle

\footnote{Zhen-Hu Ning,
 Faculty of Information Technology, Beijing University of Technology, Beijing, 100124, China. E-mail address: nzh41034@163.com.

}
\footnotetext{This work is supported by the National Science Foundation of China, grants  no.61473126 and no.61573342, and Key Research Program of Frontier Sciences, CAS, no. QYZDJ-SSW-SYS011.}
\begin{quote}
\begin{small}
{\bf Abstract} \,\,
We consider the following nonlinear Schr\"{o}dinger  equation on  exterior domain. 
\begin{equation}
\label{arma.1}\begin{cases}iu_t+\Delta_g u+ia(x)u-|u|^{p-1}u=0\qquad (x,t)\in \Om\times
(0,+\infty),
\cr u\big|_\Ga=0\qquad t\in (0,+\infty),
\cr u(x,0)=u_0(x)\qquad x\in \Om,\end{cases}
\end{equation}
where  $1<p<\frac{n+2}{n-2}$, $\Om\subset\R^n(n\ge3)$ is an exterior domain and  $(\R^n,g)$ is a complete Riemannian manifold.  We establish Morawetz  estimates  for  the system $(\ref{arma.1})$  without dissipation ($a(x)\equiv 0$ in $(\ref{arma.1})$) and  meanwhile  prove exponential stability of the system $(\ref{arma.1})$ with a dissipation effective on a neighborhood of the infinity.

It is worth mentioning that our results are different from the existing studies. First,  Morawetz  estimates  for  the system $(\ref{arma.1})$  are directly derived from the metric $g$ and are independent on the assumption of an (asymptotically) Euclidean metric. In addition, we not only prove exponential stability of the system $(\ref{arma.1})$ with non-uniform energy decay rate, which is dependent on the initial data, but also prove  exponential stability of the system $(\ref{arma.1})$ with uniform energy  decay rate.
The main methods are the development of Morawetz multipliers in non (asymptotically) Euclidean spaces  and  compactness-uniqueness arguments.
\\[3mm]
\\[3mm]
\end{small}
\end{quote}
\vskip .5cm
\tableofcontents
\def\theequation{1.\arabic{equation}}
\setcounter{equation}{0}
\section{Introduction}
\vskip .2cm
\subsection{Notations}
\quad  \ \ Let $O$ be the origin of $\R^n$ ($n\ge 3$) and
 \be r(x)= |x|,\quad  x\in\R^n\ee
be  the standard distance function of $\R^n$.
Moreover, let $\<\cdot,\cdot\>$, $\div$, $\nabla$, $\Delta$ and $I_n=(\delta_{i,j})_{n\times n}$ be the standard inner product of $\R^n$,  the standard divergence operator of $\R^n$, the standard gradient operator of $\R^n$, the standard Laplace operator of $\R^n$ and the unit matrix,respectively.

Suppose that $(\R^n,g)$ is a smooth complete Riemannian manifold with
\be g=\sum^n_{i,j=1}g_{ij}(x)dx_idx_j,\quad x\in \R^n.\ee


Let \be  G(x)=(g_{ij}(x))_{n\times n},\quad x\in \R^n.\ee
Denote
\be\<X,Y\>_g=\<G(x)X,Y\>,\quad|X|_g^2=\<X,X\>_g,\quad X,Y\in\R^n_x,x\in\R^n.\ee
Let $D$ be the Levi-Civita connection of the metric $g$ and $H$ be  a vector field,
then the covariant differential $DH$ of the vector field $H$ is a tensor field of rank 2 as follow:
\be DH(X,Y)(x)=\<D_YH,X\>_g(x)\quad  X,Y\in\R^n_x, x\in\R^n.\ee

Let $S(r)$ be the sphere in $\R^n$ with radius $r$. Then
\be\left\<X,\frac{\partial}{\partial r}\right\>=0,\quad  \textmd{for }\ \ X\in S(r)_x, x\in \R^n \backslash O.\ee

Finally, we set $\div_g$, $\nabla_g$ and $\Delta_g$  as the divergence operator of $(\R^n,g)$,  the gradient operator  of $(\R^n,g)$ and the Laplace$-$Beltrami operator of $(\R^n,g)$, respectively.

\vskip .2cm
\subsection{Nonlinear Schr\"{o}dinger  equation}

\quad \  Let $\Om\subset\R^n$ be an exterior domain with  smooth compact boundary $\Gamma$ and let $\nu(x)$ be the   unit normal vector outside  $\Om$   in $(\R^n,g)$  for  $ x\in  \Ga$. Assume that the origin $O\notin \overline{\Om}$.
Denote
\be d_1= \inf_{x\in \Ga} |x|\quad and \quad d_2=\sup_{x\in \Ga} |x|.\ee
Then $d_2\ge d_1>0$. For any constant $h>d_2$, we define
\be \Om(h)=\{x|x\in \Om,|x|\leq h\}.\ee

\quad \  We consider the following system:  
\begin{equation}
\label{wg.1} \begin{cases}iu_t+\Delta_g u+ia(x)u-|u|^{p-1}u=0\qquad (x,t)\in \Om\times
(0,+\infty),\cr u\big|_\Ga=0\qquad t\in (0,+\infty),
\cr u(x,0)=u_0(x)\qquad x\in \Om,\end{cases}
\end{equation}
where 
\be 1<p<\frac{n+2}{n-2},\ee
and $a(x)\in C^{2}(\overline{\Om})$ is a nonnegative real function satisfying
\be \label{scd.42}\sup_{x\in\overline{\Om}}\left(a(x)+\Big|\nabla_g a(x) \Big|_g +\Big|\Delta_g a(x) \Big| \right)<+\infty.\ee

\label{}

Define the energy of the system (\ref{wg.1}) by
\begin{equation}
  E(t)=\frac12\int_{\Om}\left(|u|^2+|\nabla_gu|_g^2\right)dx_g+\frac{1}{p+1}\int_{\Om}|u|^{p+1}dx_g,
\end{equation}
where
\be dx_g=\sqrt{\det(G(x))}dx,\quad |u|^2=u\bar{u},\quad |\nabla_gu|_g^2=\<\nabla_g u,\nabla_g \bar{u}\>_g.\ee

  For the free Schr\"{o}dinger  equation  on a Riemannian manifold, many Strichartz estimates  and local energy estimates are given by  \cite{BH12,Bou11, w20,1ww,ww1,w31,w3,w2,Rau78,Tsu84, 76} under the non-trapping assumption and the assumption  of  an Euclidean metric at infinity. There exists a wealth of literature on such estimates  for the wave equation (see \cite{Bur98, GHS13, w11, MRS77,Ral69, Tat13} and references therein).

 For the linear damped Schr\"{o}dinger  equation  on Riemannian manifolds,  the local energy decay in an exterior domain has been proved in \cite{Alo08b,Alo08a, AK07,AK10,Bur04,KMRJ01,Roy14} and many others under the geometric control condition (see \cite{BLR92, RT74}). Under the non-trapping condition on an exterior domain,  exponential decay for the global energy has been proved in \cite{BC14} for the Schr{\"o}dinger equation with a dissipation effective on a neighborhood of the infinity.  For the nonlinear damped Schr\"{o}dinger  equation  on compact manifold or Euclidean space,  many stability results are given by \cite{bur04-1,bur04-2,cal01,cal02,cal03,cal04,cal05,DGL01} and references therein. Such results are also based on the  non-trapping assumption or geometric control condition.

The non-trapping assumption and geometric control condition  are very closely related to the geodesic escape. Since the geodesic   depends on a nonlinear ODE, they are hard to check. On the other hand, the non-trapping assumption and geometric control condition are not sufficient to derive Morawetz estimates  for hyperbolic equations on global space. In comparison to the existing studies,  we here take advantage of the metric $g$ to establish  Morawetz estimates  for the Schr\"{o}dinger equation.

 As is known, the multiplier method is a simple and effective tool to deal with the energy  estimate on PDEs. In particular,  the celebrated Morawetz multipliers introduced by \cite{w6} have been extensively used to study the energy decay of the wave equation with
constant coefficients, see \cite{w25,w11,w16,w7} and many others. For bounded domains,  Yao\cite{w12} developed Morawetz multipliers  for the  wave equation with variable coefficients, which is a powerful tool in the analysis of systems with variable coefficients and has been extended by \cite{32,w100,30} and  many others  mentioned in \cite{w15}. However, how to establish the  Morawetz estimates in  non (asymptotically) Euclidean spaces is still an open problem. Therefore, one purpose
of this paper is to establish Morawetz estimates on  non (asymptotically) Euclidean spaces.
The organization of our paper goes as follows. In Section 2, we will state our main results. Then some multiplier identities and key lemmas for problem (\ref{wg.1}) will be presented in Section 3. We will show Morawetz estimates for the nonlinear Schr\"{o}dinger equation without dissipation  in Section 4.  Then proofs of stability of the damped nonlinear Schr\"{o}dinger  equation with non-uniform decay rate will be presented in  Section 5. We will prove stability of the damped nonlinear Schr\"{o}dinger  equation with uniform decay rate in Section 6. Finally, the proof for Assumption (U1) and Assumption (U2) hereinafter under stronger geometric condition is given in Appendix.

 \vskip .5cm
\def\theequation{2.\arabic{equation}}
\setcounter{equation}{0}
\section{Main results}

\vskip .2cm
\subsection{Well-posedness }
\quad \ \ \  Denote
\be C^{\infty}_1(\Om)=\{w\in C^{\infty} (\overline{\Om})\quad and \quad \int_{\Omega}  |w|^2 dx_g<+\infty\}.\ee
\be C^{\infty}_2(\Om)=\{w\in C^{\infty} (\overline{\Om}) \quad and \quad \int_{\Om} \left(|w|^2 +|\nabla_g w|_g^2\right) dx_g <+\infty \}.\ee
\be C^{\infty}_3(\Om)=\{w\in C^{\infty} (\overline{\Om}) \quad and \quad \int_{\Om} \left(|w|^2 +|\nabla_g w|_g^2+|\Delta_gu|^2\right) dx_g <+\infty \}.\ee
Let $L^2(\Om)$ be the closure of $C^{\infty}_1(\Om)$ with respect to the tolopogy
\be  \|w(x)\|^2_{L^2(\Om)}=\int_{\Omega} |w|^2 dx_g,\ee
 $H^1(\Om)$ be the closure of $C^{\infty}_2(\Om)$ with respect to the tolopogy
\be  \|w(x)\|^2_{H^1(\Om)}=\int_{\Om} \left(|w|^2 +|\nabla_g w|_g^2\right)dx_g,\ee
and  $H^2(\Om)$ be the closure of $C^{\infty}_3(\Om)$ with respect to the tolopogy
\be  \|w(x)\|^2_{H^2(\Om)}=\int_{\Om} \left(|w|^2 +|\nabla_g w|_g^2+|\Delta_gu|^2\right)dx_g.\ee
Denote
\be H^1_{\Ga}(\Om)=\{w\in H^1(\Om),\quad w\big|_\Ga=0\}. \ee

It is well-known that the system (\ref{wg.1}) is subcritical and  has been studied extensively in the Euclidean geometry for large classes of nonlinearities, see the books \cite{caz01,gri01}, and the references
therein. On the hyperbolic spaces, well-posedness and scattering of the system (\ref{wg.1}) without dissipation  have been proved in \cite{Ion01,Ion02}.  Therefore, throughout the paper, we assume that the following condition holds true.

\quad \  {\bf Assumption (S)}\,\,\ The  system (\ref{wg.1}) is well-posed such that
\be  u\in C\left([0,+\infty),H_{\Ga}^1(\Om)\bigcap H^2(\Om) \right).\ee

\vskip .2cm
\subsection{Morawetz estimates for the nonlinear Schr\"{o}dinger  equation  in  non (asymptotically) Euclidean spaces }
\quad \ \  The main  geometric conditions for Morawetz estimates of the nonlinear Schr\"{o}dinger  equation  in  non (asymptotically) Euclidean spaces are given by the following assumption.

{\bf Assumption (A)}\,\,\ Assume that
 \be \label{scd.23} a(x)\equiv 0\quad  in \quad \Om,\ee
 \be \label{scd.13}  G(x)\frac{\partial}{\partial r}=\frac{\partial}{\partial r}, \quad x\in \R^n,\ee
\be  \label{wg-kl.69} \left\< \left(\left(1-\alpha(x)\right)G(x)+\frac{ r}{2}\frac{\partial G(x)}{\partial r}\right)X,X \right\>\ge 0\quad for \quad X\in S(r)_x,\ \ x\in\Om,\ee
 \be  \label{wg-kl.68} \det(G(x))=c_0 r^d , \quad x\in \Om,\ee
where $c_0>0,d$ are constants and $\alpha(x)$ is a continuous nonnegative function defined on $\R^n$.
\begin{rem} Let $(r,\theta)$=$(r,\theta_1,\theta_2,\cdots,\theta_{n-1})$ be the polar coordinates of $x\in \R^n$ in the Euclidean metric. From (\ref{scd.13}), we have
 \be g=dr^2+\sum_{i,j=1}^{n-1}\gamma_{ij}(r,\theta)d\theta_id\theta_j,\quad x\in \R^n,\ee
 which implies $r(x)=|x|$ is the  geodesic distance function of $(\R^n,g)$ from $x$ to  the origin $O$.
\end{rem}
\begin{rem}\label{anapde.07} Let Assumption  ${\bf (A)}$ hold true. It follows from relations (\ref{wg-kl.74})  and   (\ref{wg-kl.75})  hereinafter that \begin{eqnarray} \frac{(n+d/2-1)}{r}&&= \frac{n-1}{r}+\frac{\partial\ln \sqrt{\det(G(x))}}{\partial r}=\Delta_g r =tr D^2r\nonumber\\
&&\ge (n-1)\frac{\alpha(x)}{r}\ge 0,\quad x\in \Om.     \end{eqnarray}
Then \be d\ge  2(1-n).\ee
\end{rem}
\begin{exl}
Let $d_1= d_2$ and  $G(x)$  satisfy
\be G(x)=\frac{x\otimes x}{|x|^2}+f(r)\left(I_n-\frac{x\otimes x}{|x|^2}\right), \quad   x\in\R^n,\ee
where $f(r)$ is a smooth function defined on $[0,+\infty)$ such that
\be f(r)=r^{m},\quad |x|\ge d_1\quad and \quad   f(r)=1,\quad |x|<\frac{d_1 }{2}.\ee
Therefore,
  \be  G(x)\frac{\partial}{\partial r} = \frac{\partial}{\partial r} ,\quad x\in \R^n,\ee
 \be  \left\< \left(\frac{1}{2}\frac{\partial G(x)}{\partial r}\right)X,X \right\>= \frac{ m}{2r} |X|^2_g\quad for \quad X\in S(r)_x,\ \ |x|\ge d_1, \ee
\be  \det(G(x))=r^{m(n-1)} \quad for \quad   |x|\ge d_1.\ee
Let \be \alpha(x)=1-\frac{ m}{2} ,\quad d=m(n-1).\ee
Then, (\ref{scd.13}),(\ref{wg-kl.69}) and (\ref{wg-kl.68}) hold true.
\end{exl}

    \begin{thm}\label{wg-kl.70} Let Assumption  ${\bf (A)}$ hold true. Assume that
       \be \label{scd.18}\frac{\partial r}{\partial \nu}\leq 0,\quad x\in \Ga.\ee
Then there exists a positive constant $C$ such that
for  $d= 2(3-n)$,
\begin{equation}
\label{wg-kl.71}
\int_0^T\int_{\Om}\frac{|u|^{p+1}}{r} dx_g dt+\int_0^T \int_{\Om} \frac{ \alpha(x)}{r} (|\nabla_{g}u|^2_g-|u_r|^2)dx_g dt \leq C E(0),
\end{equation}
and for   $d> 2(3-n)$,
\begin{equation}
\label{wg-kl.72}
 \int_0^T\int_{\Om}\left(\frac{|u|^2}{r^3} +\frac{|u|^{p+1}}{r}\right)dx_gdt+\int_0^T \int_{\Om} \frac{ \alpha(x)}{r} (|\nabla_{g}u|^2_g-|u_r|^2)dx_g dt \leq C E(0).
\end{equation}
\end{thm}

\vskip .2cm
\subsection{Stability of the damped nonlinear Schr\"{o}dinger  equation with non-uniform energy decay rate }

 \quad \   The main  geometric conditions for stability of the damped nonlinear Schr\"{o}dinger  equation with non-uniform energy decay rate  are given by the following assumption.

 {\bf Assumption (B)}\,\,\,  There exist constants  $ R_0>d_2, 0<\delta\leq 1$ such that
\be \label{2wg.7.4} \left\< \left((1-\delta)G(x)+\frac{ r}{2}\frac{\partial G(x)}{\partial r}\right)X,X \right\>\ge 0\quad for \quad X\in \R^n_x,\ \ x\in \OM(R_0),\ee
 and $a(x)$ satisfies
\be  a(x)\ge a_0>0,\quad x\in \left(\Om\backslash \Om(R_0-\varepsilon_0)\right)\bigcup \Ga(\varepsilon_1),\ee
  for some $0<2\varepsilon_1<\varepsilon_0<R_0-d_2$, where
  \be\Ga(\varepsilon)= \bigcup_{x\in \Ga} \{y\in \Om \Big| \  |y-x|<\varepsilon\},\ee
and for any $\epsilon>0$, there exists $C_\epsilon$ such that
 \be \label{wg-kl.35}\Big|\Delta_g a(x) \Big| \leq C_\epsilon a(x)+\epsilon, \quad x\in \OM.\ee

To prove the stability of the system (\ref{wg.1}), the following assumptions are also considered.

{\bf Assumption (U1)}\,\,\,   Let $\widehat{\Om}\subset \R^n$ be a  bounded domain with smooth boundary and $\omega$ be an open subset of $\widehat{\Om}$ such that
  \be \omega \supset \bigcup_{x\in \partial \widehat{\Om} } \{y\in \widehat{\Om} \Big| \  |y-x|<\xi\},\ee
  for some $\xi>0$.  Assume that $ \omega$ satisfies geometric control condition:

 {\bf (GCC)} There exists constant $T_0>0$ such that for any $x\in \widehat{\Omega}$ and any unit-speed geodesic $\gamma (t)$ of $(\R^n,g)$ starting at $x$, there exists $t< T_0$ such that $\gamma (t)\subset\omega$.

Then there exists $T_1\ge0$ such that for any $T>T_1$, the only solution $u$ in $  C([0,T],H^1(\widehat{\Om}))$  to the system
\begin{equation}
\label{shro.1}
\begin{cases}iu_t+\Delta_g u=0\qquad (x,t)\in \widehat{\Om}\times(0,T),
\cr  u=0\qquad(x,t)\in \omega\times
(0,T),\end{cases}
\end{equation}
is the trivial one $u\equiv 0$.

{\bf Assumption (U2)}\,\,\,   Let $\widehat{\Om}\subset \R^n$ be a  bounded domain with smooth boundary and $\omega$ be an open subset of $\widehat{\Om}$ such that
  \be \omega \supset \bigcup_{x\in \partial \widehat{\Om} } \{y\in \widehat{\Om} \Big| \  |y-x|<\xi\},\ee
  for some $\xi>0$.  Assume that $ \omega$ satisfies geometric control condition:

  {\bf (GCC)} There exists constant $T_0>0$ such that for any $x\in \widehat{\Omega}$ and any unit-speed geodesic $\gamma (t)$ of $(\R^n,g)$ starting at $x$, there exists $t< T_0$ such that $\gamma (t)\in \omega$.

 Then there exists $T_1\ge 0$ such that for any $T>T_1$, the only  solution $u$ in $C([0,T],H^1(\widehat{\Om}))$  to the system
\begin{equation}
\label{shro.2}
\begin{cases}iu_t+\Delta_g u-|u|^{p-1}u=0\qquad (x,t)\in \widehat{\Om}\times(0,T),
\cr  u=0\qquad(x,t)\in \omega\times
(0,T).\end{cases}
\end{equation}
is the trivial one $u\equiv 0$.

\begin{rem}
 If $T_1=0$, which implies $T$ can be arbitrary small in (\ref{shro.1}) and  (\ref{shro.2}), Assumption  ${\bf (U1)}$  and  Assumption  ${\bf (U2)}$ are called as  unique continuation condition. On Euclidean space, unique continuation condition for  linear(or nonlinear) Schr\"{o}dinger  equation has been proved by \cite{Escauriaza10, Escauriaza12,Ion06, Kenig102, Koch05, Wang1} and the references therein. On  Riemannian
manifold,  under the assumption that  unique continuation condition for  linear Schr\"{o}dinger  equation holds true,   unique continuation condition for  the nonlinear Schr\"{o}dinger  equation was  proved by \cite{Laurent110} in  dimension $3$ and \cite{DGL01} in dimension $2$ .

By the equivalent relation between the controllability and the observability estimate \cite{Lions1988},
Assumption  ${\bf (U1)}$ follows from  Theorem 4.4 in  \cite{Laurent14}. However, a detailed proof of Theorem 4.4 in \cite{Laurent14} is not provided.

 Under a stronger geometric condition than (GCC), we can prove Assumption  ${\bf (U1)}$  and  Assumption  ${\bf (U2)}$ directly by multiplier methods. See  Proposition (\ref{Appendix.1}) and  Proposition (\ref{Appendix.5}) in the Appendix.
\end{rem}

    \begin{thm}\label{wg1.7.2_1} Let Assumption  ${\bf (B)}$,  Assumption  ${\bf (U1)}$ and  Assumption  ${\bf (U2)}$ hold true. Assume that  $\Big|\Big|u_0\Big|\Big|_{L^2(\Om)}  \leq E_0.$
Then there exist  positive constants $C_1$ and $C_2$, which are dependent on $E_0$, such that
\be E(t)\leq C_1e^{-C_2t} E(0),\quad \forall t>0.\ee\end{thm}

\vskip .2cm
\subsection{Stability of the damped nonlinear Schr\"{o}dinger  equation  with uniform energy decay rate }

\quad \ The main  geometric conditions for stability of the damped nonlinear Schr\"{o}dinger  equation with uniform energy decay rate  are given by the following assumption.

 {\bf Assumption (C)}\,\,\,There exist constants  $R_0>d_2, 0<\delta\leq 1$  such that
\be \label{wg-kl.63}  G(x)\frac{\partial}{\partial r}=\frac{\partial}{\partial r} ,\quad |x|\leq R_0 \quad and \quad \det(G(x))= c_0 r^d ,\quad  x\in \OM(R_0),\ee
\be \label{wg-kl.64}  \left\< \left((1-\delta)G(x)+\frac{ r}{2}\frac{\partial G(x)}{\partial r}\right)X,X \right\>\ge 0\quad for \quad X\in \R^n_x,\ \ x\in \OM(R_0),\ee
where $c_0>0,d$ are constants and $a(x)$ satisfies
    \be  a(x)\ge a_0>0,\quad x\in\Om\backslash \Om(R_0-\varepsilon_0),\ee
 for some $0<\varepsilon_0<R_0-d_2$ and for any $\epsilon>0$, there exists $C_\epsilon$ such that
 \be \label{wg-kl.38}\Big|\Delta_g a(x) \Big| \leq C_\epsilon a(x)+\epsilon, \quad x\in \OM.\ee
\begin{rem} Let Assumption  ${\bf (C)}$ hold true. It follows from the relations(\ref{wg-kl.74})  and   (\ref{wg-kl.75}) hereinafter that \begin{eqnarray} \frac{(n+d/2-1)}{r}&&= \frac{n-1}{r}+\frac{\partial\ln \sqrt{\det(G(x))}}{\partial r}=\Delta_g r =tr D^2r\nonumber\\
&&\ge (n-1)\frac{\delta}{r},\quad x\in \Om(R_0).     \end{eqnarray}
Then \be d\ge  2(n-1)(\delta-1).\ee
\end{rem}

   \begin{thm}\label{wg-kl.37} Let Assumption  ${\bf (C)}$ hold true. Assume that
       \be\label{scd.20} \frac{\partial r}{\partial \nu}\leq 0,\quad x\in \Ga.\ee
Then there exist  positive constants $C_1,C_2$ such that
\be E(t)\leq C_1e^{-C_2t} E(0),\quad \forall t>0.\ee  \end{thm}

 \vskip .5cm
\def\theequation{3.\arabic{equation}}
\setcounter{equation}{0}
\section{Multiplier Identities and Key Lemmas}
\vskip .2cm

\quad \ \ We need to establish several multiplier identities, which are useful for our problem.

\begin{lem}\label{wg.14}
Let $\widehat{\Om}\subset\R^n$ be a bounded domain with smooth boundary. Suppose that $u(x,t)$ solves the following equation:
\be \label{wg-kl.92}iu_t+\Delta_g u+ia(x)u-|u|^{p-1}u=0\qquad (x,t)\in \widehat{\Om} \times
(0,+\infty).\ee
Let $\H$ be a  $C^{1}$vector
field defined on $\overline{\widehat{\Om}}$.  Then
\begin{eqnarray}
 \label{wg.14.1}
 &&\int_0^T\int_{\partial \widehat{\Om}}
\Re \left(\frac{\pa u}{\pa\hat{\nu}}\H(\bar{u}) \right)d\Ga_g dt+\frac12\int_0^T\int_{\partial \widehat{\Om}}
\left(\Im (u\bar{u}_t)-\left|\nabla_g u\right|_g^2-\frac{2}{p+1}|u|^{p+1}\right)\<\H,\hat{\nu}\>_g d\Ga_g dt\nonumber\\
=&&\frac{1}{2}\int_{\widehat{\Om}}\Im\left(u \H(\bar{u})\right) dx_g\Big |^T_0+\int_0^T\int_{\widehat{\Om}}\Re D\H(\nabla_g \bar{u},\nabla_g u)  dx_g dt\nonumber\\
\quad &&+\int_0^T\int_{\widehat{\Om}} \Im\left(a(x) u\H(\bar{u}) \right) dx_g dt\nonumber\\
\quad &&+\frac12\int_0^T\int_{\widehat{\Om}}\left(\Im (u\bar{u}_t)-\left|\nabla_g u\right|_g^2-\frac{2}{p+1}|u|^{p+1}\right)\div_g\H dx_g dt,
\end{eqnarray}
where $\hat{\nu}(x)$  is the   unit normal vector outside  $\widehat{\Om}$   in $(\R^n,g)$  for  $ x\in  \partial \widehat{\Om}$ and
\be d\Ga_g=\sqrt{\det(G(x))}d\Ga.\ee

Moreover, assume that the real function $P\in C^2(\overline{\widehat{\Om}})$. Then
\begin{eqnarray}
\label{wg.14.2}
&&\int_0^T\int_{\widehat{\Om}}\left(\Im (u\bar{u}_t)-\left|\nabla_g u\right|_g^2-|u|^{p+1}\right)P dx_g dt\nonumber\\
&& = \frac12\int_0^T\int_{\partial \widehat{\Om}}|u|^2\frac{\pa
P}{\pa\hat{\nu}}d\Ga_g dt-\frac12 \int_0^T
\int_{\widehat{\Om}}|u|^2(\Delta_g P) dx_g dt\nonumber\\
&&\qquad-\int_0^T\int_{\partial \widehat{\Om}}\Re(P \bar{u}\frac{\pa
u}{\pa\hat{\nu}})d\Ga_g dt.\end{eqnarray}
\end{lem}

{\bf Proof}.
  Firstly, we multiply  (\ref{wg-kl.92}) by $\H(\bar{u})$ and integrate over $\widehat{\Om}\times
(0,T) $. We deduce that
\begin{eqnarray} \Re\left(iu_t \H(\bar{u})\right)=&& -\Im \left(u_t \H(\bar{u})\right)\nonumber\\
=&& -\frac{1}{2}\Im \left(u_t \H(\bar{u})-\bar{u}_t \H(u)\right)\nonumber\\
=&&-\frac{1}{2}\Im \left((u \H(\bar{u}))_t- \H(u\bar{u}_t)\right)\nonumber\\
=&&-\frac{1}{2}\Im \left(u \H(\bar{u})\right)_t+\frac{1}{2}\Im \H(u\bar{u}_t)\nonumber\\
=&&-\frac{1}{2}\Im \left(u \H(\bar{u})\right)_t+\frac{1}{2}\Im \div_g(u\bar{u}_t\H) -\frac{1}{2}\Im(u\bar{u}_t\div_g\H),
\end{eqnarray}
\begin{eqnarray} \Re \left(\H(\bar{u})\Delta_g u)\right)=&&\Re \left( \div_g \H(\bar{u}) \nabla_g u-\nabla_g u\<\H,\nabla_g \bar{u}\>_g \right)\nonumber\\
=&&  \Re  \div_g \H(\bar{u}) \nabla_g u- \Re \nabla_g u\<\H,\nabla_g \bar{u}\>_g   \nonumber\\
=&&  \Re  \div_g \H(\bar{u}) \nabla_g u- \Re D\H(\nabla_g \bar{u},\nabla_g u) - \Re D^2\bar{u}(\H,\nabla_g u) \nonumber\\
=&& \Re  \div_g \H(\bar{u}) \nabla_g u- \Re D\H(\nabla_g \bar{u},\nabla_g u) - \Re D^2\bar{u}(\nabla_g u,\H) \nonumber\\
=&&  \Re  \div_g \H(\bar{u}) \nabla_g u- \Re D\H(\nabla_g \bar{u},\nabla_g u) - \frac{1}{2}\H(|\nabla_g u|^2_g)\nonumber\\
=&&  \Re  \div_g \H(\bar{u}) \nabla_g u- \Re D\H(\nabla_g \bar{u},\nabla_g u) - \frac{1}{2}\div_g (|\nabla_g u|^2_g\H)\nonumber\\
&&\quad +\frac{1}{2} |\nabla_g u|^2_g\div_g\H,
\end{eqnarray}
and
\begin{eqnarray} \Re\left(ia(x)u-|u|^{p-1}u\right)\H(\overline{u})
&&=-\Im\left(a(x)u\H(\overline{u})\right)\nonumber\\
 &&\quad- \frac{1}{p+1} \div_g \left(|u|^{p+1} \H\right)+\frac{|u|^{p+1}}{p+1} \div_g \H.
\end{eqnarray}

The equality  (\ref{wg.14.1}) follows from Green's formula.

In addition, by multiplying (\ref{wg-kl.92})  by $P\bar{u}$ and integrating over $\widehat{\Om}\times
(0,T)$, we obtain
\begin{eqnarray} \Re \left(i P u_t \bar{u}\right)=&& -\Im \left(P u_t \bar{u}\right)= \Im \left(P u \bar{u}_t \right),
\end{eqnarray}
\begin{eqnarray} \Re \left(P\bar{u}\Delta_g u\right)=&&\Re \left( \div_g P\bar{u} \nabla_g u-\nabla_g u(P\bar{u}) \right)\nonumber\\
=&&  \Re  \div_g P\bar{u} \nabla_g u-P|\nabla_gu|^2_g -\frac{1}{2}\nabla_g P(|u|^2)  \nonumber\\
=&&  \Re \div_g P\bar{u} \nabla_g u-P|\nabla_gu|^2_g -\frac{1}{2}\div_g |u|^2\nabla_g P+\frac{1}{2}|u|^2\Delta_g P ,
\end{eqnarray}
and
\begin{eqnarray} \Re\left(ia(x)u-|u|^{p-1}u\right)P\overline{u}&&
=\Re\left(ia(x)P|u|^2\right)- P|u|^{p+1}
 \nonumber\\ &&=- P|u|^{p+1}.
\end{eqnarray}

 The equality  (\ref{wg.14.2}) follows from Green's formula.$\Box$

The following  lemma will be utilized frequently in our subsequent proof.
 \begin{lem}\label{pde.4}  Let $x_0\in\R^n$ be a fixed point. Let $H(x)=x-x_0$, then
 \be \label{pde.5} DH(X,X)=\left\<\left(G(x)+\frac{ \widehat{r}(x)}{2}\frac{\partial G(x)}{\partial \widehat{r}}\right)X,X\right\>,\quad \textmd{for}\ \  X\in \R^n_x, x \in \R^n,\ee
 where $\widehat{r}(x)=|x-x_0|$.
  \end{lem}

  {\bf Proof}. Let $x\in\R^n,  X=\sum_{i=1}^nX_i\frac{\partial}{\partial x_i}\in\R^n_x$.
  Note that
  \be  H(x)=\sum_{i=1}^n(x_i-x_{0,i})\frac{\partial}{\partial x_i}.\ee
  Then, we deduce that
 \begin{eqnarray} DH(X,X)=&&\sum_{i,j,k=1}^n\left\< D_{\frac{\partial}{\partial x_i}}\left((x_k-x_{0,k})\frac{\partial}{\partial x_k}\right),\frac{\partial}{\partial x_j}\right\>_{g}X_iX_j \nonumber\\
=&& \sum_{i,j=1}^n g_{ij}X_iX_j+\sum_{i,j,k=1}^n (x_k-x_{0,k})\left\< D_{\frac{\partial}{\partial x_i}}\frac{\partial}{\partial x_k},\frac{\partial}{\partial x_j}\right\>_gX_iX_j\nonumber\\
=&& |X|^2_g+\sum_{i,j,k=1}^n(x_k-x_{0,k})\left\< D_{\frac{\partial}{\partial x_k}}\frac{\partial}{\partial x_i},\frac{\partial}{\partial x_j}\right\>_gX_iX_j\nonumber\\
=&& |X|^2_g+\sum_{i,j,k=1}^n\frac{ (x_k-x_{0,k})}{2}\frac{\partial g_{ij}}{\partial x_k}X_iX_j\nonumber\\
=&&\left\<\left(G(x)+\frac{ \widehat{r}(x)}{2}\frac{\partial G(x)}{\partial \widehat{r}}\right)X,X\right\>.\end{eqnarray}

$\Box$

The following lemmas show the  relationship between the metric $g$ and  geometric control condition.
\begin{lem}\label{wg-kl.58}  Let $\widehat{\Omega}\subset\R^n$be a bounded domain and  $x_0\in\R^n$ be a fixed point. Assume that there exists $\delta>0$  such that
 \be \left\<\left((1-\delta)G(x)+\frac{ \widehat{r}(x)}{2}\frac{\partial G(x)}{\partial \widehat{r}}\right)X,X\right\>\ge 0,\quad \textmd{for}\ \  X\in \R^n_x, x \in \overline{\widehat{\Omega}},\ee
 where  $\widehat{r}(x)=|x-x_0|$.
Then,  for any $x\in \widehat{\Omega}$ and any unit-speed geodesic $\gamma (t)$ starting at $x$,
if
   \be \gamma (t)\in\widehat{\Omega},\quad 0\leq t \leq t_0, \ee
then
\be t_0 \leq \frac{2}{\delta}\sup\left\{|x-x_0|_g(x)\Big|\ \  x\in \overline{\widehat{\Omega}}\right\}.\ee
\end{lem}

{\bf Proof.} Let $H(x)=x-x_0$. It follows from (\ref{pde.5})  that
\be D H(X,X)\geq
\delta|X|_g^2\quad for\ \ all\ \ X\in\R^n_x,\quad
x\in\overline{\widehat{\Omega}}.\ee

Note that
\be |\gamma' (t)|_g=1,\quad D_{\gamma' (t)}\gamma' (t)=0. \ee
Then
\be\<H,\gamma' (t)\>_g\Big|_0^{t_0}= \int_0^{t_0}\gamma' (t)\<H,\gamma' (t)\>_gdt=\int_0^{t_0}D H(\gamma' (t),\gamma' (t))dt\ge \delta t_0. \ee
Hence
  \be t_0 \leq\frac{2}{\delta}\sup\left\{|H|_g(x)\Big|\ \  x\in \overline{\widehat{\Omega}}\right\}.\ee $\Box$

\begin{lem}\label{cmp.3_4}
Assume that
 \be  G(x)\frac{\partial}{\partial r}=\frac{\partial}{\partial r}, \quad x\in \R^n,\ee
    \be \frac{\partial G(x)}{\partial r}=- \frac{2}{r_2}G(x)\left(I_n-\frac{x\otimes x}{|x|^2}\right), \quad |x|=r_2,\ee
    where $r_2$ is a positive constant.
Then, for any $x\in S(r_2)$ and any unit-speed geodesic $\gamma (t)$ starting at $x$ with
 \be\gamma' (0)\in S(r_2)_x,\ee
 we have
  \be\gamma(t)\in S(r_2),\quad \forall t\ge 0.\ee

\end{lem}
{\bf Proof }\,\,\,
Note that
\be D(rDr)=Dr\otimes Dr+rD^2r.\ee
With  (\ref{pde.5}),  we obtain
\be D^2r(X,X)=\left\< \left(\frac{1}{r}G(x)+\frac{ 1}{2}\frac{\partial G(x)}{\partial r}\right)X,X \right\>=0\quad for \quad X\in S(r_2)_x,\ \ x\in S(r_2).\ee

    Let $\widehat{g}$ be a Riemannian metric induced by $g$ in $S(r_2)$ and $\widehat{D}$ be the associated Levi-Civita connection.

    Let $\widehat{\gamma}(t)$ be a unit-speed geodesic of  $(S(r_2),\widehat{g})$ starting at $x\in S(r_2)$, then
 \be \left\<{\widehat{\gamma}'(t)},\frac{\partial}{\partial r}\right\>_g=0,\quad \widehat{D}_{\widehat{\gamma}'(t)}{\widehat{\gamma}'(t)}=0,\quad \forall t\ge0.\ee
Therefore,
 \begin{eqnarray}  D_{\widehat{\gamma}'(t)}{\widehat{\gamma}'(t)}
&&=\widehat{D}_{\widehat{\gamma}'(t)}{\widehat{\gamma}'(t)}+\left\<D_{\widehat{\gamma}'(t)}{\widehat{\gamma}'(t)},\frac{\partial}{\partial r}\right\>_g\frac{\partial}{\partial r}
\nonumber\\
&& =\widehat{D}_{\widehat{\gamma}'(t)}{\widehat{\gamma}'(t)}-D^2r(\widehat{\gamma}'(t), \widehat{\gamma}'(t))\frac{\partial}{\partial r}
=0, \end{eqnarray}
which implies  $\widehat{\gamma}(t)$ is also a geodesic of $(\R^n,g)$.
$\Box$

\vskip .5cm
\def\theequation{4.\arabic{equation}}
\setcounter{equation}{0}
\section{Proofs of Morawetz estimates in  non (asymptotically) Euclidean spaces}

\begin{lem}\label{wg-kl.22}Let $u(x,t)$ solve the system (\ref{wg.1}). Then
\be \label{wg-kl.23} \int_{\Om }|u|^2dx_g \Big|_0^T=-2\int_0^T \int_{\Om  }a(x) |u|^2dx_g dt ,\ee
\begin{eqnarray}\label{wg-kl.45} \int_{\Om  }\left(|\nabla_gu|^2_g+ \frac{2}{p+1}|u|^{p+1} \right) dx_g \Big|_0^T&&=-2\int_0^T \int_{\Om  }a(x)\left( |\nabla_gu|^2_g+|u|^{p+1} \right)dx_g dt\nonumber\\
&& + \int_0^T
\int_{\Om}|u|^2(\Delta_g a(x)) dx_g dt,\end{eqnarray} for any $T>0$.\end{lem}
{\bf Proof}.  Multiplying the Schr\"{o}dinger equation in (\ref{wg.1}) by $2\bar{u}$ and integrating over $\Om\times
(0,T)$, we have
\be \int_{\Om  }|u|^2dx_g \Big|_0^T=-2\int_0^T \int_{\Om  }a(x) |u|^2dx_g .\ee

After multiplying the Schr\"{o}dinger equation in (\ref{wg.1}) by $2\bar{u}_t$ and then integrating over $\Om\times
(0,T)$, we obtain
\be \label{wg2.20.5.1}  \int_{\Om  }\left(|\nabla_gu|^2_g+ \frac{2}{p+1}|u|^{p+1} \right) dx_g \Big|_0^T =-2 \int_0^T \int_{\Om  } \Im(a(x)u \bar{u}_t)dx_gdt.\ee
Let $P=a(x)$ and $\widehat{\Om}=\Om(a)$ in (\ref{wg.14.2}). Substituting (\ref{wg.14.2})  into (\ref{wg2.20.5.1}), letting $a\rightarrow +\infty$, we get
\begin{eqnarray}\int_{\Om  }\left(|\nabla_gu|^2_g+ \frac{2}{p+1}|u|^{p+1} \right) dx_g \Big|_0^T&&=-2\int_0^T \int_{\Om  }a(x)\left( |\nabla_gu|^2_g+|u|^{p+1} \right)dx_g dt\nonumber\\
&& + \int_0^T
\int_{\Om}|u|^2(\Delta_g a(x)) dx_g dt.\end{eqnarray}  $\Box$

\begin{lem}\label{wg-kl.73} Let Assumption  ${\bf (A)}$ hold true. Then
\be \label{wg-kl.74} D^2r(X,X)\geq
\frac{\alpha(x)}{r}|X|_g^2\quad \textmd{for}\ \ \textmd{all}\ \  X\in S(r)_x,
x\in\Om ,\ee
\be \label{wg-kl.75}\Delta_g r=\frac{n+d/2-1}{r}\quad \textmd{for} \ \ x\in \Om .\ee
  \end{lem}
{\bf Proof.} Note that
\be D(rDr)=Dr\otimes Dr+rD^2r.\ee
With  (\ref{pde.5}),  we obtain
\be D^2r(X,X)=\left\< \left(\frac{1}{r}G(x)+\frac{ 1}{2}\frac{\partial G(x)}{\partial r}\right)X,X \right\>\geq
\frac{\alpha(x)}{r}|X|_g^2\quad for \quad X\in S(r)_x,\ \ x\in\Om ,\ee
and
\be\Delta_g r= \frac{n-1}{r}+\frac{\partial\ln \sqrt{\det(G(x))}}{\partial r}=\frac{n+d/2-1}{r}\quad \textmd{for} \ \ x\in \Om .\ee $\Box$\\

{\bf Proof  of Theorem \ref{wg-kl.70}}

$\quad \ $Let $\H=\frac{\partial}{\partial r}$ and $\widehat{\Om}=\Om(h)$ in (\ref{wg.14.1}).
It follows from  (\ref{wg.14.1}), (\ref{wg-kl.74}) and  (\ref{wg-kl.75}),  that
\begin{eqnarray}
 \label{wg-kl.76}
\int_0^T&&\int_{\partial \Om(h)}
\Re  \left(\frac{\pa u}{\pa\hat{\nu}}\H(\bar{u}) \right)d\Ga_g dt \nonumber\\
&&\quad+\frac12\int_0^T\int_{\partial \Om(h)}
\left(\Im(u\bar{u}_t)-\left|\nabla_g u\right|_g^2-\frac{2}{p+1}|u|^{p+1}\right)\<\H,\hat{\nu}\>_g d\Ga_g dt \nonumber\\
 \geq&&\frac{1}{2}\int_{\Om(h)}\Im \left(u \bar{u}_r\right) dx_g\Big |^T_0+\int_0^T\int_{\Om(h)}\frac{\alpha(x)}{r}(|\nabla_{g}u|^2_g-|u_r|^2)dx_g dt\nonumber\\
&&\quad+\int_0^T\int_{\Om(h)}\left(\Im (u\bar{u}_t)-\left|\nabla_g u\right|_g^2-\frac{2}{p+1}|u|^{p+1}\right)\frac{n-1}{2r} dx_g dt \nonumber\\
 =&&\frac{1}{2}\int_{\Om(h)}\Im \left(u \bar{u}_r\right) dx_g\Big |^T_0+\int_0^T\int_{\Om(h)}\frac{\alpha(x)}{r}(|\nabla_{g}u|^2_g-|u_r|^2)dx_g dt\nonumber\\
&&\quad +\int_0^T\int_{\Om(h)}\left(\Im (u\bar{u}_t)-\left|\nabla_g u\right|_g^2-|u|^{p+1}\right)\frac{n-1}{2r} dx_g dt\nonumber\\
&& \quad+\int_0^T\int_{\Om(h)}\frac{(n-1)(p-1)}{2r(p+1)}|u|^{p+1} dx_g dt.
\end{eqnarray}
Let $P=\frac{n-1}{2r}$ and $\widehat{\Om}=\Om(h)$ in (\ref{wg.14.2}). Substituting (\ref{wg.14.2})  into (\ref{wg-kl.76}), letting $h\rightarrow +\infty$, we obtain
\begin{eqnarray}
 \label{wg-kl.77}\frac12\int_{\Om}&&\Im \left(u  \bar{u}_r\right) dx_g\Big |^T_0- \frac{n-1}{4}\int_0^T
\int_{\Om}|u|^2 \Delta_g\left(\frac{1 }{r} \right) dx_g dt\nonumber\\
 &&\quad +\int_0^T\int_{\Om} \frac{(p-1)(n-1)}{2(p+1)r}|u|^{p+1} dx_g dt+\int_0^T\int_{\Om}\frac{\alpha(x)}{r}(|\nabla_{g}u|^2_g-|u_r|^2)dx_g dt\nonumber\\
 &&\le \Pi_{\Ga},
\end{eqnarray}
where
\begin{eqnarray}
 \label{scd.15} \Pi_{\Ga}&&=\int_0^T\int_{\Ga}
\Re  \left(\frac{\pa u}{\pa\nu}\H(\bar{u}) \right)d\Ga_g dt \nonumber\\
&&\quad +\frac12\int_0^T\int_{\Ga}
\left(\Im(u\bar{u}_t)-\left|\nabla_g u\right|_g^2-\frac{2}{p+1}|u|^{p+1}\right)\<\H,\nu\>_g d\Ga_g dt \nonumber\\
&&\qquad-\frac12\int_0^T\int_{ \Ga}|u|^2\frac{\pa
P}{\pa\nu}d\Ga_g dt+\int_0^T\int_{ \Ga}\Re(P \bar{u}\frac{\pa
u}{\pa\nu})d\Ga_g dt.
\end{eqnarray}

 Since $u\Large|_{\Ga}=0,$ we
obtain $\nabla_{\Ga_g} \bar{u}\Big|_{\Ga}=0$, that is,
\be\label{scd.17} \nabla_g \bar{u}=\frac{\pa
 \bar{u}}{\pa\nu}\nu\quad\mbox{for}\quad
 x\in\Ga.\ee
Similarly, we have \be \label{scd.16} \H(\bar{u})=\<\H,\nabla_g \bar{u}\>_g=\frac{\pa
 \bar{u}}{\pa\nu}\<\H,\nu\>_g=\frac{\pa
 \bar{u}}{\pa\nu}\frac{\pa
 r}{\pa\nu}\quad\mbox{for}\quad
 x\in\Ga.\ee
Using the formulas (\ref{scd.17}) and   (\ref{scd.16}) in the formula
(\ref{scd.15}) on the portion $\Ga$, with (\ref{scd.18}), we obtain
\begin{equation}
\label{scd.21}\Pi_{\Ga}=\frac12\int_0^T
\int_{\Ga}\Big|\frac{\pa
u}{\pa\nu}\Big|^2\frac{\pa
 r}{\pa\nu} d\Ga_g dt\le 0.
\end{equation}
Substituting (\ref{scd.21})  into (\ref{wg-kl.77}), we have
\begin{eqnarray}
\label{scd.22}\frac12\int_{\Om}&&\Im \left(u  \bar{u}_r\right) dx_g\Big |^T_0- \frac{n-1}{4}\int_0^T
\int_{\Om}|u|^2 \Delta_g\left(\frac{1 }{r} \right) dx_g dt\nonumber\\
 &&\quad +\int_0^T\int_{\Om} \frac{(p-1)(n-1)}{2(p+1)r}|u|^{p+1} dx_g dt+\int_0^T\int_{\Om}\frac{\alpha(x)}{r}(|\nabla_{g}u|^2_g-|u_r|^2)dx_g dt\nonumber\\
 &&\le 0.
\end{eqnarray}

Note that
\be  \Delta_g \left(\frac{1}{r}\right) =-\frac{n+d/2-3}{r}.\ee
With  (\ref{scd.22}), we obtain
for  $d= 2(3-n)$,
\begin{eqnarray}
 \label{wg-kl.78}\frac{1}{2}\int_{\Om}&&\Im \left(u  \bar{u}_r\right) dx_g\Big |^T_0 +\int_0^T\int_{\Om} \frac{(p-1)(n-1)}{2(p+1)r}|u|^{p+1} dx_g dt \nonumber\\
 &&\quad+\int_0^T\int_{\Om}\frac{\alpha(x)}{r}(|\nabla_{g}u|^2_g-|u_r|^2)dx_g dt\nonumber\\
 &&\le 0,
\end{eqnarray}
and for  $d> 2(3-n)$,
\begin{eqnarray}
 \label{wg-kl.79}
\frac{1}{2}\int_{\Om}&&\Im \left(u  \bar{u}_r\right) dx_g\Big |^T_0+\int_0^T\int_{\Om}\frac{(n-1)(n+d/2-3)}{4 r^3}u^2 dx dt\nonumber\\
 &&\quad +\int_0^T\int_{\Om} \frac{(p-1)(n-1)}{2(p+1)r}|u|^{p+1} dx_g dt+\int_0^T\int_{\Om}\frac{\alpha(x)}{r}(|\nabla_{g}u|^2_g-|u_r|^2) dx dt\nonumber\\
 &&\le 0.
\end{eqnarray}

It follows from (\ref{scd.23}), (\ref{wg-kl.23}) and (\ref{wg-kl.45}) that
\be E(t)= E(0),\quad t>0.\ee
The estimates (\ref{wg-kl.71}) and  (\ref{wg-kl.72}) follows from (\ref{wg-kl.78}) and (\ref{wg-kl.79}).
$\Box$\\

\vskip .5cm
\def\theequation{5.\arabic{equation}}
\setcounter{equation}{0}
\section{Proofs of stability  with non-uniform energy decay rate}

\quad \ From Lemma \ref{wg-kl.58}, the following lemma holds true.
\begin{lem} Let assumption  ${\bf (B)}$ hold true.
Then, there exists $t_0>0$,  for any $x\in  \OM(R_0)$ and any unit-speed geodesic $\gamma (t)$ starting at $x$, there exists $t<t_0$ such that
\be \label{scd.01} \gamma (t)\in  \partial \OM(R_0).\ee
\end{lem}

\begin{lem}\label{wg-kl.28}Let assumption  ${\bf (B)}$ hold true. Let $u(x,t)$  solve the system (\ref{wg.1}). Then
\begin{eqnarray}\label{wg-kl.27}
E(0)+\int_0^TE(t) dt&&\leq C\int_0^T \int_{\Om  }a(x) \left(|u|^2+|\nabla_gu|^2_g+|u|^{p+1}\right)dx_g dt\nonumber\\
&&\quad+C\int_0^T\int_{\Om(R_0-\varepsilon_0)} |u|^2 dx_g dt,\end{eqnarray}  for sufficiently large $T$.\end{lem}
{\bf Proof}.
Let $b(x)\in C^{\infty}(\R^n)$  be a nonnegative  function  satisfying
\be b(x)=1,\quad x\in \Om(R_0-\varepsilon_0)\backslash  \Ga(\varepsilon_1)      \quad and \quad b(x)=0, \quad  x\in \R^n\backslash \Om(R_0).\ee

Let  \be \label{scd.41}H(x)= b(x)x,\quad  x\in \R^n.\ee
It follows from  (\ref{2wg.7.4}) and (\ref{pde.5}) that
\be\label{wg1.20.1} D H(X,X)\geq
\delta |X|_g^2 \quad \textmd{for}\ \ \textmd{all}\ \  X\in \R^n_x,
x\in \Om(R_0-\varepsilon_0)\backslash  \Ga(\varepsilon_1),\ee
\be\div_g H=tr DH \ge n\delta
 \quad \textmd{for}\ \ \textmd{all}\ \
x\in \Om(R_0-\varepsilon_0)\backslash  \Ga(\varepsilon_1).\ee

Let $\H= H $ and $\widehat{\Om}=\OM(R_0)$ in (\ref{wg.14.1}). From (\ref{wg.14.1}), we have
\begin{eqnarray}
 \label{wg.20.3}
  0 \geq&&\frac{1}{2}\int_{\OM(R_0)}\Im \left(u H(\bar{u})\right) dx_g\Big |^T_0+\delta\int_0^T\int_{\Om(R_0-\varepsilon_0)\backslash  \Ga(\varepsilon_1)}|\nabla_{g}u|^2_gdx_g dt\nonumber\\
&&\quad -C\int_0^T\int_{ \left(\Om(R_0)\backslash \Om(R_0-\varepsilon_0)\right)\bigcup   \Ga(\varepsilon_1)}|\nabla_{g}u|^2_gdx_g dt+\int_0^T\int_{\OM(R_0)} \Im\left(a(x) uH(\bar{u}) \right) dx_g dt\nonumber\\
&&\quad +\frac12\int_0^T\int_{\OM(R_0)}\left(\Im (u\bar{u}_t)-\left|\nabla_g u\right|_g^2- \frac{2}{p+1}|u|^{p+1}  \right)\div_g H dx_g dt\nonumber\\
=&&  \frac{1}{2}\int_{\OM(R_0)}\Im \left(u H(\bar{u})\right) dx_g\Big |^T_0+\delta\int_0^T\int_{\Om(R_0-\varepsilon_0)\backslash  \Ga(\varepsilon_1)}|\nabla_{g}u|^2_gdx_g dt\nonumber\\
&&\quad -C\int_0^T\int_{ \left(\Om(R_0)\backslash \Om(R_0-\varepsilon_0)\right)\bigcup   \Ga(\varepsilon_1)}|\nabla_{g}u|^2_gdx_g dt+\int_0^T\int_{\OM(R_0)} \Im\left(a(x) uH(\bar{u}) \right) dx_g dt\nonumber\\
&&\quad+\frac12\int_0^T\int_{\OM(R_0)}\left(\Im (u\bar{u}_t)-\left|\nabla_g u\right|_g^2- |u|^{p+1}\right)\div_g H dx_g dt\nonumber\\
&&\quad+\int_0^T\int_{\OM(R_0)}\frac{(p-1)\div_g H}{2(p+1)}|u|^{p+1} dx_g dt.
\end{eqnarray}
Let $P=\frac{\div_g H }{2}$ and $\widehat{\Om}=\OM(R_0)$ in (\ref{wg.14.2}). Substituting (\ref{wg.14.2})  into (\ref{wg.20.3}), we obtain
\begin{eqnarray}
 \label{scd.28}&&\frac{1}{2}\int_{\OM(R_0)}\Im \left(u H(\bar{u})\right) dx_g\Big |^T_0-\frac14 \int_0^T
\int_{\OM(R_0)}|u|^2 \Delta_g(\div_g H)  dx_g dt\nonumber\\
 &&\quad+\int_0^T\int_{\Om(R_0)} \Im\left(a(x) uH(\bar{u}) \right) dx_g dt\nonumber\\
 &&\quad +\delta\int_0^T\int_{\Om(R_0-\varepsilon_0) \backslash\Ga(\varepsilon_1)} |\nabla_{g}u|^2_gdx_g dt+\int_0^T\int_{\Om(R_0-\varepsilon_0)\backslash\Ga(\varepsilon_1)} \frac{n\delta(p-1)}{2(p+1)}|u|^{p+1}dx_g dt \nonumber\\
 &&\le C\int_0^T\int_{ \left(\Om(R_0)\backslash \Om(R_0-\varepsilon_0)\right)\bigcup   \Ga(\varepsilon_1)} \left( |\nabla_{g}u|^2_g+|u|^{p+1}\right)dx_g dt.
\end{eqnarray}
Therefore
\begin{eqnarray}
&&\int_0^T\int_{\Om(R_0-\varepsilon_0)\backslash\Ga(\varepsilon_1)}  \left(|\nabla_{g}u|^2_g +|u|^{p+1}\right)dx_g dt\nonumber\\
 \leq&& C(E(0)+E(T)) +C\int_0^T \int_{\OM(R_0)  }a(x) \left(|u|^2+|\nabla_gu|^2_g+|u|^{p+1}\right)dx_g dt\nonumber\\
 &&+C \int_0^T
\int_{\Om(R_0-\varepsilon_0)}|u|^2  dx_g dt.
\end{eqnarray}
Hence
\begin{eqnarray}
\label{wg-kl.24}\int_0^TE(t) dt &&\leq C(E(0)+E(T))\nonumber\\
 &&\quad +C\int_0^T \int_{\Om  }a(x) \left(|u|^2+|\nabla_gu|^2_g+|u|^{p+1}\right)dx_g dt\nonumber\\
 &&\quad +C \int_0^T
\int_{\Om(R_0-\varepsilon_0)}|u|^2  dx_g dt.
\end{eqnarray}

With (\ref{wg-kl.23}) and (\ref{wg-kl.45}), we deduce that
\begin{eqnarray}\label{wg-kl.46}
CE(T)&&\leq CE(0)+C\int_0^{T} \int_{\Om  }a(x) \left(|u|^2+|\nabla_gu|^2_g+|u|^{p+1}\right)dx_g dt\nonumber\\
&&\quad  +
\frac{C}{2} \int_0^T
\int_{\Om}|u|^2\Big|\Delta_g a(x)\Big| dx_g dt,
\end{eqnarray}
and
\begin{eqnarray}\label{wg-kl.25}
4CE(0)&&=\int_0^{4C}E(t)dt-\int_0^{4C}(E(t)-E(0))dt\nonumber\\
&&\leq \int_0^{4C}E(t)dt+4C\int_0^{4C} \int_{\Om  }a(x) \left(|u|^2+|\nabla_gu|^2_g+|u|^{p+1}\right)dx_g dt\nonumber\\
&&\quad  +
2C \int_0^T
\int_{\Om}|u|^2\Big|\Delta_g a(x)\Big| dx_g dt.
\end{eqnarray}
Substituting (\ref{wg-kl.46}) and (\ref{wg-kl.25})  into (\ref{wg-kl.24}), for $T>4C$, with (\ref{wg-kl.35}),  we have
\begin{eqnarray}\label{wg-kl.26}
E(0)+\int_0^TE(t)  dt\leq &&C\int_0^T \int_{\Om  }a(x) \left(|u|^2+|\nabla_gu|^2_g+|u|^{p+1}\right)dx_g dt\nonumber\\
&&\quad +C\int_0^T\int_{\Om(R_0-\varepsilon_0)} |u|^2 dx_g dt.
\end{eqnarray}

The estimate (\ref{wg-kl.27}) holds true.

\begin{lem}\label{wg-kl.29} Let assumption  ${\bf (B)}$,  assumption  ${\bf (U1)}$ and  assumption  ${\bf (U2)}$ hold true and  let $T$ be sufficiently large. Then for any $\Big|\Big|u_0\Big|\Big|_{L^2(\Omega)}  \leq E_0$, there exists positive constant $C(E_0,T)$ such that
\begin{eqnarray}\label{wg-kl.30}
E(0)+\int_0^TE(t) dt&&\leq C(E_0,T)\int_0^T \int_{\Om  }a(x) \left(|u|^2+|\nabla_gu|^2_g+|u|^{p+1}\right)dx_g dt.\end{eqnarray}  \end{lem}
{\bf Proof}. We apply  compactness-uniqueness arguments to prove the conclusion. It follows from (\ref{wg-kl.27}) that
 \begin{eqnarray}  \label{wg-kl.85}
E(0)+\int_0^TE(t) dt&&\leq C\int_0^T \int_{\Om  }a(x) \left(|u|^2+|\nabla_gu|^2_g+|u|^{p+1}\right)dx_g dt\nonumber\\
&&\quad+C\int_0^T\int_{\Om(R_0-\varepsilon_0)} |u|^2 dx_g dt.\end{eqnarray}
Then, if the estimate (\ref{wg-kl.30}) doesn't hold true, there exist $\Big\{u_k\Big\}_{k=1}^{\infty}$ such that
  \be \label{wg-kl.81} \int_0^T\int_{\Om(R_0-\varepsilon_0)} |u_k|^2 dx_g dt\ge k \int_0^T \int_{\Om  }a(x) \left(|u_k|^2+|\nabla_gu_k|^2_g  +|u_k|^{p+1}\right)dx_g dt.\ee
Thus,
 \be \label{scd.08} E_k(0)+\int_0^TE_k(t) dt\leq C E_0,\ee
 where \be
  E_k(t)=\frac12\int_{\Om}\left(|u_k|^2+|\nabla_gu_k|_g^2\right)dx_g+\frac{1}{p+1}\int_{\Om}|u_k|^{p+1}dx_g.\ee

 Therefore, there exists $\hat{u}_0$ and a subset of $\Big\{u_k\Big\}_{k=1}^{\infty}$, still denoted by $\Big\{u_k\Big\}_{k=1}^{\infty}$,  such that
   \be \label{anapde.01} u_k\rightarrow \hat{u}_0 \ \ weakly \ \  in  \ \   L^2([0,T],H_\Ga^1(\Om)),\ee
   and
 \be u_k\rightarrow \hat{u}_0 \ \ strongly \ \ in \ \ L^2(\Om(R_0-\varepsilon_0))\ for \ arbitrarily\ fixed\ t\in [0,T].\ee
Note that \be \|u_k-\hat{u}\|^2_{L^2(\Om(R_0-\varepsilon_0))}\le \widehat{C}(T)E_0,\ \ \forall t \in[0,T],\ \forall 1\le k< +\infty.\ee
Lebesgue's dominated convergence theorem yields
\be u_k\rightarrow \hat{u}_0 \ \ strongly \ \ in  \ \  L^2(\Om(R_0-\varepsilon_0)\times(0,T)).\ee
{\bf Case a:}
  \be \label{wg-kl.80} \int_0^T\int_{\Om(R_0-\varepsilon_0)} |\hat{u}_0|^2 dx_g dt>0.\ee

  It follows from  (\ref{scd.42}), (\ref{wg-kl.23}), (\ref{wg-kl.45}) and (\ref{scd.08}) that there exists $C(T)>0$ such that
\be \label{anapde.06} E_k(t)\leq C(T)E_0 ,\quad \forall 0\leq t\leq T.\ee

Denote
\be \label{anapde.02} q=\frac{2n}{(n-2)p},\quad q^{*}=\frac{q}{q-1}.\ee
Since $1<p<\frac{n+2}{n-2}$, then
\be \frac{2n}{n+2}<q,q^{*}<\frac{2n}{n-2}.\ee
Note that
\be\frac{1}{q}+\frac{1}{q^{*}}=1,\ee
then, $L^{q^{*}}\left(\Om(R_0-\varepsilon_0)\right)$ is the dual space of $L^q\left(\Om(R_0-\varepsilon_0)\right)$.

Note that
\be \label{scd.07}H_\Ga^1\left(\Om(R_0-\varepsilon_0)\right)\hookrightarrow L^{\frac{2n}{n-2}}\left(\Om(R_0-\varepsilon_0)\right).\ee
therefore, it follows from (\ref{anapde.06}) that
\be \{|u_k|^{p-1}u_k\}\ are \ bounded \ in\  L^{\infty}([0,T],L^{q}(\Om(R_0-\varepsilon_0))).\ee
Then
\be \{|u_k|^{p-1}u_k\}\ are \ bounded \ in\   L^q \left(\Om(R_0-\varepsilon_0)\times(0,T)\right).\ee
Hence, there exists a subset of $\Big\{u_k\Big\}_{k=1}^{\infty}$, still denoted by $\Big\{u_k\Big\}_{k=1}^{\infty}$,  such that
  \be \label{scd.43} |u_k|^{p-1}u_k \rightarrow |\hat{u}_0|^{p-1}\hat{u}_0\ \ weakly \ \  in  \ \  L^q \left(\Om(R_0-\varepsilon_0)\times(0,T)\right).\ee

 It follows from (\ref{wg-kl.81}) that
 \be a(x)\hat{u}_0=0\qquad(x,t)\in \Om\times
(0,T).\ee
 Therefore, with (\ref{anapde.01}) and  (\ref{scd.43}),  we obtain
   \begin{equation}
\label{wg-kl.32}\begin{cases}i\hat{u}_{0t}+\Delta_g \hat{u}_0-|\hat{u}_0|^{p-1}\hat{u}_0=0\qquad (x,t)\in \left(\Om(R_0-\varepsilon_0)\times(0,T)\right),
\cr a(x)\hat{u}_0=0\qquad(x,t)\in \Om\times
(0,T).\end{cases}
\end{equation}
With (\ref{scd.01}) and Assumption {\bf  (U2)}, we have
\be \label{wg-kl.82} \hat{u}_0\equiv0,\qquad (x,t)\in \Om\times
(0,T), \ee
which contradicts (\ref{wg-kl.80}).

{\bf Case b:}
 \be \label{scd.05} \hat{u}_0\equiv 0 \quad on \quad  \Om(R_0-\varepsilon_0)\times(0,T). \ee

Denote
\be \label{scd.06}v_k=  u_k\Big/ \sqrt{c_k} \quad for\quad  k\ge 1,\ee
where
\be \label{anapde.04} c_k=\int_0^T\int_{\Om(R_0-\varepsilon_0)}  |u_k|^2 dx_gdt.\ee

Then $v_k$ satisfies
   \begin{equation}
 \label{wg-kl.90}\begin{cases}iv_{kt}+\Delta_g v_k+ia(x)v_k-|u_k|^{p-1}v_k=0\qquad (x,t)\in \Om\times
(0,T),\cr v_k\Big|_\Ga=0\qquad t\in (0,T),
\end{cases}
\end{equation}
and
\be  \label{wg-kl.89}\int_0^T\int_{\Om(R_0-\varepsilon_0)} |v_k|^2 dx_g dt=1.\ee
It follows from (\ref{wg-kl.81}) that
\be \label{wg-kl.88}1\ge k \int_0^T \int_{\Om  }a(x) \left(|v_k|^2+|\nabla_gv_k|^2_g  +|u_{k}|^{p-1}|v_k|^{2}\right)dx_g dt.\ee
Therefore, it follows from  (\ref{wg-kl.85}) that
\begin{eqnarray}  \label{wg-kl.86}
\widehat{E}_k(0)+\int_0^T\widehat{E}_k(t) dt&&\leq 1+\frac{1}{k}\leq 2,\end{eqnarray}
where
 \be\widehat{E}_k(t)=\int_{\Om  } \left(|v_k|^2+|\nabla_gv_k|^2_g  +|u_{k}|^{p-1}|v_k|^{2}\right)dx_g.\ee
 Hence, there exists $v_0$ and a subset of $\Big\{v_k\Big\}_{k=1}^{\infty}$, still denoted by $\Big\{v_k\Big\}_{k=1}^{\infty}$, such that
   \be \label{anapde.05} v_k\rightarrow v_0 \ \ weakly \ \  in  \ \   L^2([0,T],H_\Ga^1(\Om)),\ee
   and
    \be v_k\rightarrow \hat{v}_0 \ \ strongly \ \ in \ \ L^2(\Om(R_0-\varepsilon_0))\ for \ arbitrarily\ fixed\ t\in [0,T].\ee
    Then by Lebesgue's dominated convergence theorem, we obtain
    \be\label{anapde1.05} v_k\rightarrow v_0 \ \ strongly \ \ in  \ \  L^2(\Om(R_0-\varepsilon_0)\times(0,T)).\ee

It follows from  (\ref{scd.42}), (\ref{wg-kl.23}) and (\ref{wg-kl.45})  that there exists $C(T)>0$ such that
\be E_k(t)\leq C(T)E_k(0) ,\quad \forall 0\leq t\leq T.\ee
With (\ref{scd.06})  and (\ref{wg-kl.86}), we obtain
\be \label{scd.44}\widehat{E}_k(t)\leq C(T) \widehat{E}_k(0)\leq 2C(T),\quad \forall 0\leq t\leq T.\ee

Let $q,q^{*}$ be given by (\ref{anapde.02}).  Note that
 \be \label{scd.07}H_\Ga^1\left(\Om(R_0-\varepsilon_0)\right)\hookrightarrow L^{\frac{2n}{n-2}}\left(\Om(R_0-\varepsilon_0)\right).\ee
 Therefore, it follows from (\ref{scd.44}) that
\be \{|v_k|^{p-1}v_k\}\ are \ bounded \ in\  L^{\infty}([0,T],L^{q}(\Om(R_0-\varepsilon_0))).\ee
Hence
\begin{eqnarray} \int_0^T\int_{\Om(R_0-\varepsilon_0)}&&\left(|u_k|^{p-1}|v_k|\right)^q dx_gdt\nonumber\\ &&=c_k^{\frac{q(p-1)}{2}} \int_0^T\int_{\Om(R_0-\varepsilon_0)} |v_k|^{\frac{2n}{n-2}} dx_gdt\nonumber\\ && \leq c_k^{\frac{q(p-1)}{2}} C(T).\end{eqnarray}
With  (\ref{scd.05}), (\ref{anapde.04})  and  (\ref{anapde1.05}), we obtain
  \be \label{anapde.03}\lim_{k\rightarrow +\infty}\int_0^T\int_{\Om(R_0-\varepsilon_0)}\left(|u_k|^{p-1}|v_k| \right)^qdx_gdt= 0.\ee

It follows from (\ref{wg-kl.88}) that
 \be a(x)v_0=0\qquad(x,t)\in \Omega\times
(0,T).\ee
Therefore, it follows from (\ref{wg-kl.90}), (\ref{anapde.05}) and (\ref{anapde.03}) that
   \begin{equation}
\begin{cases}iv_{0t}+\Delta_g v_0=0\qquad (x,t)\in \Om(R_0-\varepsilon_0)\times(0,T),
\cr  a(x)v_0=0\qquad(x,t)\in \Om\times
(0,T).\end{cases}
\end{equation}
With (\ref{scd.01}) and Assumption {\bf  (U1)}, we have
\be \label{wg-kl.91} v_0\equiv0,\qquad (x,t)\in \Om\times
(0,T). \ee

   It follows from (\ref{wg-kl.89})  that
 \be \label{wg-kl.31}\int_0^T\int_{\Om(R_0-\varepsilon_0)} |v_0|^2 dx_g dt=1,\ee
which contradicts (\ref{wg-kl.91}). $\Box$ \\

{\bf Proof of Theorem \ref{wg1.7.2_1} }\,\,\,
 Let $T$ be sufficiently large.
It follows from (\ref{wg-kl.23}) that $\Big|\Big|u\Big|\Big|_{L^2(\Om)}$ is non increasing.
Hence, with (\ref{wg-kl.30}), we obtain
\begin{eqnarray}
\label{scd.14}E(S)+&&\int_S^{S+T}E(t) dt\nonumber\\&&\quad\leq C(E_0,T)\int_S^{S+T} \int_{\Om }a(x) \left(|u|^2+|\nabla_gu|^2_g+|u|^{p+1}\right)dx_g dt,\end{eqnarray}
for any $S\ge 0$.

 It follows from  (\ref{wg-kl.23}) and (\ref{wg-kl.45}) that
 \be \int_{\Om }|u|^2dx_g \Big|_S^{S+T}=-2\int_S^{S+T} \int_{\Om  }a(x) |u|^2dx_g dt,\ee
 and
\begin{eqnarray}\int_{\Om  }\left(|\nabla_gu|^2_g+ \frac{2}{p+1}|u|^{p+1} \right) dx_g \Big|_S^{S+T}&&=-2\int_S^{S+T}\int_{\Om  }a(x)\left( |\nabla_gu|^2_g+|u|^{p+1} \right)dx_g dt\nonumber\\
&& + \int_S^{S+T}
\int_{\Om}|u|^2(\Delta_g a(x)) dx_g dt.\end{eqnarray}
 Therefore, with (\ref{scd.14}), we deduce that
  \begin{eqnarray}\label{wg-kl.33}
E(S)+&& \int_S^{S+T}E(t) dt\nonumber\\ &&\leq C(E_0,T)(E(S)-E(S+T)) +C(E_0,T)\int_S^{S+T}\int_{\Om} \Big|\Delta_g a(x) \Big||u|^2 dx_g dt\nonumber\\ &&= C(E_0,T)(E(S)-E(S+T)) -M\int_{\Om} |u|^2 dx_g\Big |_S^{S+T}\nonumber\\ &&\quad +\int_S^{S+T}\int_{\Om} \left(C(E_0,T)\Big|\Delta_g a(x) \Big|  -Ma(x)\right) |u|^2 dx_g dt.
\end{eqnarray}
For sufficiently large $M$, with (\ref{wg-kl.35})  we have
 \begin{eqnarray}\label{wg-kl.34}
E(S)\leq &&C(E_0,T)(E(S)-E(S+T)) -M\int_{\Om} |u|^2 dx_g \Big |_S^{S+T}.
\end{eqnarray}

Denote
\be \label{wg-kl.61}\widetilde{E}(t)= E(t)+\frac{M}{C(E_0,T)}\int_{\Om} |u|^2 dx_g.\ee
From  (\ref{wg-kl.34}),  we obtain
\begin{eqnarray}\label{wg-kl.36}
 &&\widetilde{E}(S)\leq\widetilde{C}(E_0,T)(\widetilde{E}(S)-\widetilde{E}(S+T)).
\end{eqnarray}
Then
\be \label{wg-kl.62} \widetilde{E}(S+T) \leq  \frac{\widetilde{C}(E_0,T)-1}{\widetilde{C}(E_0,T)} \widetilde{E}(S).\ee

It follows from (\ref{scd.42}), (\ref{wg-kl.23}), (\ref{wg-kl.45}) and (\ref{wg-kl.61})  that there exists $\widetilde{C}(T)>0$ such that
\be \widetilde{E}(S+t)\leq \widetilde{C}(T)  \widetilde{E}(S),\quad \forall 0\leq t\leq T.\ee
With (\ref{wg-kl.62}), $\widetilde{E}(t)$ is of exponential decay.
Hence, there exist $C_1(E_0),C_2(E_0)>0$ such that
\be  E(t) \leq C_1(E_0)e^{-C_2(E_0) t} E(0), \forall t>0.\ee $\Box$

\vskip .5cm
\def\theequation{6.\arabic{equation}}
\setcounter{equation}{0}
\section{Proofs of stability  with uniform energy decay rate}

\begin{lem}\label{wg-kl.39} Let Assumption  ${\bf (C)}$ hold true. Assume that
       \be\label{scd.20} \frac{\partial r}{\partial \nu}\leq 0,\quad x\in \Ga.\ee Let $u(x,t)$  solve the system (\ref{wg.1}). Then
\begin{eqnarray}\label{wg-kl.40}
E(0)+\int_0^TE(t) dt\leq &&C\int_0^T \int_{\Om  }a(x) \left(|u|^2+|\nabla_gu|^2_g+|u|^{p+1}\right)dx_g dt,\end{eqnarray}  for sufficiently large $T$.\end{lem}
{\bf Proof}.\ \ \
Let $b(z)$  be a smooth nonnegative  function  defined on $[0,+\infty)$  satisfying
\be \label{wg-kl.84} b(z)=1,\quad  0\leq z\leq R_0-\varepsilon_0 \quad and \quad b(z)=0,\quad  z\ge  R_0.  \ee
Let
 \be\label{wg-k1l.47} H(x)= b(r)x,\quad x\in \R^n.\ee
It follows from  (\ref{wg-kl.63}), (\ref{wg-kl.64}) and (\ref{pde.5}) that
\be\label{wg-kl.47} D H(X,X)\geq
\delta|X|_g^2  \quad \textmd{for}\ \ \textmd{all}\ \  X\in \R^n_x,
x\in \Om(R_0-\varepsilon_0),\ee
\begin{eqnarray} \div_gH &&=1+r\Delta_gr\nonumber\\
&&=1+r\left(\frac{n-1}{r}+\frac{\partial\ln \sqrt{\det(G(x))}}{\partial r}\right)\nonumber\\
&&=n+d/2 \quad \textmd{for}\ \ \textmd{all}\ \  X\in \R^n_x,
x\in \Om(R_0-\varepsilon_0),\end{eqnarray}

Let $\H= H $ and $\widehat{\Om}=\OM(R_0)$ in (\ref{wg.14.1}). From (\ref{wg.14.1}), we have
\begin{eqnarray}
 \label{wg-kl.48}
\int_0^T&&\int_{\Ga}
\Re  \left(\frac{\pa u}{\pa\nu}\H(\bar{u}) \right)d\Ga_g dt \nonumber\\
&&\quad +\frac12\int_0^T\int_{\Ga}
\left(\Im(u\bar{u}_t)-\left|\nabla_g u\right|_g^2-\frac{2}{p+1}|u|^{p+1}\right)\<\H,\nu\>_g d\Ga_g dt
\nonumber\\ \ge&& \frac{1}{2}\int_{\OM(R_0)}\Im \left(u H(\bar{u})\right) dx_g\Big |^T_0+\delta\int_0^T\int_{\Om(R_0-\varepsilon_0)}|\nabla_{g}u|^2_gdx_g dt\nonumber\\
&&\quad-C\int_0^T\int_{x\in \OM(R_0)\backslash \Om(R_0-\varepsilon_0)}|\nabla_{g}u|^2_gdx_g dt+\int_0^T\int_{\OM(R_0)} \Im\left(a(x) uH(\bar{u}) \right) dx_g dt\nonumber\\
&&\quad+\frac12\int_0^T\int_{\OM(R_0)}\left(\Im (u\bar{u}_t)-\left|\nabla_g u\right|_g^2- \frac{2}{p+1}|u|^{p+1}\right)\div_g H dx_g dt\nonumber\\ =&&\frac{1}{2}\int_{\OM(R_0)}\Im \left(u H(\bar{u})\right) dx_g\Big |^T_0+\delta\int_0^T\int_{\Om(R_0-\varepsilon_0)}|\nabla_{g}u|^2_gdx_g dt\nonumber\\
&&\quad-C\int_0^T\int_{x\in \OM(R_0)\backslash \Om(R_0-\varepsilon_0)}|\nabla_{g}u|^2_gdx_g dt+\int_0^T\int_{\OM(R_0)} \Im\left(a(x) uH(\bar{u}) \right) dx_g dt\nonumber\\
&&\quad+\frac12\int_0^T\int_{\OM(R_0)}\left(\Im (u\bar{u}_t)-\left|\nabla_g u\right|_g^2- |u|^{p+1}\right)\div_g H dx_g dt\nonumber\\
&&\quad+\int_0^T\int_{\OM(R_0)}\frac{(p-1)\div_g H}{2(p+1)}|u|^{p+1} dx_g dt.
\end{eqnarray}
Let $P=\frac{\div_g H }{2}$ and $\widehat{\Om}=\OM(R_0)$  in (\ref{wg.14.2}). Substituting (\ref{wg.14.2})  into (\ref{wg-kl.48}), we obtain
\begin{eqnarray}
 \label{scd.34}\frac{1}{2}\int_{\OM(R_0)}&&\Im \left(u H(\bar{u})\right) dx_g\Big |^T_0-\frac14 \int_0^T
\int_{\OM(R_0)}|u|^2 \Delta_g(\div_g H)  dx_g dt\nonumber\\
 &&\quad+\int_0^T\int_{\OM(R_0)} \Im\left(a(x) u H(\bar{u}) \right) dx_g dt +\delta\int_0^T\int_{\Om(R_0-\varepsilon_0)} |\nabla_{g}u|^2_gdx_g dt\nonumber\\
 &&\quad+\int_0^T\int_{\Om(R_0-\varepsilon_0)} \frac{(n+d/2)(p-1)}{2(p+1)}|u|^{p+1}dx_g dt \nonumber\\
\le&& C\int_0^T\int_{\OM(R_0)\backslash \Om(R_0-\varepsilon_0)} \left(|\nabla_gu|^2_g+|u|^{p+1}\right)dx_g dt+ \Pi_{\Ga},
\end{eqnarray}
where
\begin{eqnarray}
 \label{scd.30} \Pi_{\Ga}=&&\int_0^T\int_{\Ga}
\Re  \left(\frac{\pa u}{\pa\nu}H(\bar{u}) \right)d\Ga_g dt \nonumber\\
&&\quad +\frac12\int_0^T\int_{\Ga}
\left(\Im(u\bar{u}_t)-\left|\nabla_g u\right|_g^2-\frac{2}{p+1}|u|^{p+1}\right)\<H,\nu\>_g d\Ga_g dt \nonumber\\
&&\quad-\frac12\int_0^T\int_{ \Ga}|u|^2\frac{\pa
P}{\pa\nu}d\Ga_g dt+\int_0^T\int_{ \Ga}\Re(P \bar{u}\frac{\pa
u}{\pa\nu})d\Ga_g dt.
\end{eqnarray}

 Since $u\Large|_{\Ga}=0,$ we
obtain $\nabla_{\Ga_g} \bar{u}\Big|_{\Ga}=0$, that is,
\be\label{scd.31} \nabla_g \bar{u}=\frac{\pa
 \bar{u}}{\pa\nu} \nu\quad\mbox{for}\quad
 x\in\Ga.\ee
Then, with (\ref{wg-k1l.47}), we have \be \label{scd.32} H(\bar{u})=\<H,\nabla_g \bar{u}\>_g=\frac{\pa
 \bar{u}}{\pa\nu}\<H,\nu\>_g=r\frac{\pa
 \bar{u}}{\pa\nu}\frac{\pa
 r}{\pa\nu}\quad\mbox{for}\quad
 x\in\Ga.\ee
Using the formulas (\ref{scd.31}) and   (\ref{scd.32}) in formula
(\ref{scd.30}) on the portion $\Ga$, with (\ref{scd.20}), we obtain
\begin{equation}
\label{scd.33}\Pi_{\Ga}=\frac12\int_0^T
\int_{\Ga}r\Big|\frac{\pa
u}{\pa\nu}\Big|^2\frac{\pa
 r}{\pa\nu} d\Ga_g dt\le 0.
\end{equation}
Substituting (\ref{scd.33})  into (\ref{scd.34}), we have
\begin{eqnarray}
 \label{scd.35}&&\frac{1}{2}\int_{\OM(R_0)}\Im \left(u H(\bar{u})\right) dx_g\Big |^T_0-\frac14 \int_0^T
\int_{\OM(R_0)}|u|^2 \Delta_g(\div_g H)  dx_g dt\nonumber\\
 &&\quad+\int_0^T\int_{\OM(R_0)} \Im\left(a(x) uH(\bar{u}) \right) dx_g dt\nonumber\\
 &&\quad +\delta\int_0^T\int_{\Om(R_0-\varepsilon_0)} |\nabla_{g}u|^2_gdx_g dt+\int_0^T\int_{\Om(R_0-\varepsilon_0)} \frac{(n+d/2ㄘ(p-1)}{(p+1)}|u|^{p+1}dx_g dt \nonumber\\
 &&\le C\int_0^T\int_{\OM(R_0)\backslash \Om(R_0-\varepsilon_0)} \left( |\nabla_{g}u|^2_g+|u|^{p+1}\right)dx_g dt.
\end{eqnarray}
Therefore,
\begin{eqnarray}
\label{wg-kl.50}\int_0^T\int_{\Om(R_0-\varepsilon_0)}&& \left( |u|^{p+1}+|\nabla_{g}u|^2_g\right)dx_g dt\leq C(E(0)+E(T))\nonumber\\
 &&\quad +C\int_0^T \int_{\OM(R_0) }a(x) \left(|u|^2+|\nabla_gu|^2_g+|u|^{p+1}\right)dx_g dt.
\end{eqnarray}

Note that
\be \int_0^T\int_{\Om(R_0-\varepsilon_0)} |u|^{2}dx_g dt\leq C(R_0)\int_0^T\int_{\Om(R_0-\varepsilon_0)} |\nabla_gu|^2_gdx_g dt.\ee
Hence
\begin{eqnarray}
\label{wg-kl.56}\int_0^T E(t) dt &&\leq C(E(0)+E(T))\nonumber\\
 &&\quad +C\int_0^T \int_{\Om  }a(x)\left(|u|^2+|\nabla_gu|^2_g+|u|^{p+1}\right)dx_g dt.
\end{eqnarray}

With (\ref{wg-kl.23}) and  (\ref{wg-kl.45}) , we deduce that
\begin{eqnarray}
\label{scd.09}
CE(T)&&\leq CE(0)+C\int_0^{T} \int_{\Om  }a(x) \left(|u|^2+|\nabla_gu|^2_g+|u|^{p+1}\right)dx_g dt\nonumber\\
&&\quad  +
\frac{C}{2} \int_0^T
\int_{\Om}|u|^2\Big|\Delta_g a(x)\Big| dx_g dt,
\end{eqnarray}
and
\begin{eqnarray}
\label{scd.10}
4CE(0)&&=\int_0^{4C}E(t)dt-\int_0^{4C}(E(t)-E(0))dt\nonumber\\
&&\leq \int_0^{4C}E(t)dt+4C\int_0^{4C} \int_{\Om  }a(x) \left(|u|^2+|\nabla_gu|^2_g+|u|^{p+1}\right)dx_g dt\nonumber\\
&&\quad  +
2C \int_0^T
\int_{\Om}|u|^2\Big|\Delta_g a(x)\Big| dx_g dt.
\end{eqnarray}
Substituting (\ref{scd.09}) and (\ref{scd.10})  into (\ref{wg-kl.56}), for $T>4C$, with (\ref{wg-kl.38}),  we have
\begin{eqnarray}\label{wg-kl.51}
E(0)+\int_0^TE(t) dt\leq &&C\int_0^T \int_{\Om  }a(x)\left(|u|^2+|\nabla_gu|^2_g+|u|^{p+1}\right)dx_g dt.
\end{eqnarray}

 The estimate (\ref{wg-kl.40}) holds true. $\Box$\\

{\bf Proof of Theorem \ref{wg-kl.37} }\,\,\,
From (\ref{wg-kl.23}), (\ref{wg-kl.45}) and (\ref{wg-kl.40}), we deduce that
 \begin{eqnarray}
E(0)+\int_0^T&& E(t) dt\leq C(E(0)-E(T)) +C\int_0^T\int_{\Om(R_0-\varepsilon_0)} \Big|\Delta_g a(x) \Big||u|^2 dx_g dt\nonumber\\ &&= C(E(0)-E(T)) -M\int_0^T\int_{\Om} |u|^2 dx_g dt\Big |^T_0\nonumber\\ &&\quad +\int_0^T\int_{\Om(R_0-\varepsilon_0)} \left(C\Big|\Delta_g a(x) \Big|  -M\int_0^T\int_{\Om}a(x)\right) |u|^2 dx_g dt.
\end{eqnarray}
For sufficiently large $M$, with (\ref{wg-kl.38})  we have
 \begin{eqnarray}\label{wg-kl.57}
E(0)\leq &&C(E(0)-E(T)) -M\int_0^T\int_{\Om} |u|^2 dx_g dt\Big |^T_0.
\end{eqnarray}

Denote
\be \label{wg-kl.65}\widetilde{E}(t)= E(t)+M\int_0^T\int_{\Om} |u|^2 dx_g dt.\ee
From  (\ref{wg-kl.57}),  we obtain
\begin{eqnarray}\label{wg-kl.66}
 &&\widetilde{E}(0)\leq\widetilde{C}(\widetilde{E}(0)-\widetilde{E}(T)).
\end{eqnarray}
Then
\be \label{wg-kl.67} \widetilde{E}(T) \leq  \frac{\widetilde{C}-1}{\widetilde{C}} \widetilde{E}(0).\ee

It follows from  (\ref{scd.42}), (\ref{wg-kl.23}), (\ref{wg-kl.45}) and (\ref{wg-kl.65})  that there exists $\widetilde{C}(T)>0$ such that
\be \widetilde{E}(t)\leq \widetilde{C}(T)  \widetilde{E}(0),\quad \forall 0\leq t\leq T.\ee
With (\ref{wg-kl.67}), $\widetilde{E}(t)$ is exponentially decaying.
Hence, there exist $C_1,C_2>0$ such that
\be  E(t) \leq C_1e^{-C_2 t} E(0), \forall t>0.\ee $\Box$

\subsection*{Acknowledgements }
The authors would like to express their gratitude to the editors and reviewers for
their valuable comments and helpful suggestions.

  \appendix
  \renewcommand{\appendixname}{Appendix~\Alph{section}}
   \section{Appendix: Proofs of  Assumption (U1) and Assumption (U2) under a strong geometric condition }
  \quad \ \  Let $\widehat{\Om}\subset \R^n$ be a bounded domain with smooth boundary and $\omega$ be an open subset of $\widehat{\Om}$ such that
  \be \omega \supset \bigcup_{x\in \partial \widehat{\Om} } \{y\in \widehat{\Om} \Big| \  |y-x|<\xi\},\ee
  for some $\xi>0$.

   Assume that  the origin $O\notin  \overline{\widehat{\Om}}$ and
\be \label{Appendix.2} G(x)\frac{\partial}{\partial r}=\frac{\partial}{\partial r} , \quad x\in \R^n, \quad and \quad \det(G(x))= c_0 r^d ,\quad  x\in  \widehat{\Om},\ee
\be \label{Appendix.3}  \left\< \left((1-\delta)G(x)+\frac{ r}{2}\frac{\partial G(x)}{\partial r}\right)X,X \right\>\ge 0\quad for \quad X\in \R^n_x,\ \ x\in \widehat{\Om},\ee
where  $0<\delta\leq 1$, $c_0>0$ and $d$ are constants.
\begin{rem} It follows from  (\ref{wg-kl.74})  and   (\ref{wg-kl.75})  that \begin{eqnarray} \frac{(n+d/2-1)}{r}&&= \frac{n-1}{r}+\frac{\partial\ln \sqrt{\det(G(x))}}{\partial r}=\Delta_g r =tr D^2r\nonumber\\
&&\ge (n-1)\frac{\delta}{r},\quad x\in \widehat{\Om}.     \end{eqnarray}
Then \be d\ge  2(n-1)(\delta-1).\ee
\end{rem}

   \begin{pro}\label{Appendix.1} There exists $T_1\ge0$ such that for any $T>T_1$, the only solution $u$ in $  C([0,T],H^1(\widehat{\Om}))$  to the system
\begin{equation}
\label {Appendix.4}
\begin{cases}iu_t+\Delta_g u-|u|^{p-1}u=0\qquad (x,t)\in \widehat{\Om}\times(0,T),
\cr  u=0\qquad(x,t)\in \omega\times
(0,T),\end{cases}
\end{equation}
is the trivial one $u\equiv 0$. \end{pro}
 {\bf Proof}.\ \ \
Let $b(x)\in C^{\infty}(\R^n)$  be a nonnegative  function satisfying
\be  b(x)=1,\quad  x\in \widehat{\Om}\backslash \omega  \quad and \quad b(x)=0,\quad \R^n \backslash \widehat{\Om}.  \ee
Let
 \be H(x)= b(x)x,\quad x\in \R^n.\ee
It follows from  (\ref{Appendix.2}), (\ref{Appendix.3}) and (\ref{pde.5}) that
\be  D H(X,X)\geq
\delta|X|_g^2  \quad \textmd{for}\ \ \textmd{all}\ \  X\in \R^n_x,
x\in  \widehat{\Om}\backslash \omega,\ee
\begin{eqnarray} \div_gH &&=1+r\Delta_gr\nonumber\\
&&=1+r\left(\frac{n-1}{r}+\frac{\partial\ln \sqrt{\det(G(x))}}{\partial r}\right)\nonumber\\
&&=n+d/2 \quad \textmd{for}\ \ \textmd{all}\ \  X\in \R^n_x,\
x\in  \widehat{\Om}\backslash \omega.\end{eqnarray}

Let $a(x)=0$ in (\ref{wg-kl.92}).
Let $\H= H $ in (\ref{wg.14.1}) and  $P=\frac{\div_g H }{2}$ in (\ref{wg.14.2}).  Substituting (\ref{wg.14.2})  into (\ref{wg.14.1}), we obtain
\begin{eqnarray}
\frac{1}{2}\int_{\widehat{\Om}}&&\Im \left(u H(\bar{u})\right) dx_g\Big |^T_0-\frac14 \int_0^T
\int_{\widehat{\Om}}|u|^2 \Delta_g(\div_g H)  dx_g dt\nonumber\\
 &&\quad+\int_0^T\int_{\widehat{\Om}} \Re\left(DH(\nabla_g \bar{u},\nabla_g u) \right)  dx_g dt\nonumber\\
 &&\quad+\int_0^T\int_{\widehat{\Om}} \frac{(p-1)\div_gH}{2(p+1)}|u|^{p+1}dx_g dt \nonumber\\
&&=0.
\end{eqnarray}
Then
\begin{eqnarray}
&&\int_0^T\int_{\widehat{\Om}\backslash\omega}\left(|\nabla_g u|^2+\frac{2}{p+1}|u|^{p+1}\right)  dx_g dt \nonumber\\
\leq&&  C\Big|\int_{\widehat{\Om}}\Im \left(u H(\bar{u})\right) dx_g\Big |^T_0\Big|
+C\int_0^T\int_{\omega}\left(|\nabla_g u|^2+|u|^2+|u|^{p+1}\right)  dx_g dt.
\end{eqnarray}
Hence
\begin{eqnarray}
\int_0^T\int_{\widehat{\Om}}\left(|\nabla_g u|^2+\frac{2}{p+1}|u|^{p+1}\right)  dx_g dt&&\le2CE(0) .
\end{eqnarray}
Note that
 \be \int_0^T\int_{\widehat{\Om}} |u|^{2}dx_g dt\leq C\int_0^T \int_{\widehat{\Om}} |\nabla_gu|^2_g dx_g dt.\ee
Therefore
\be  \int_0^T\int_{\widehat{\Om}}\left(|u|^2+|\nabla_g u|^2+\frac{2}{p+1}|u|^{p+1}\right)  dx_g dt \leq 2CE(0),  \ee
which implies
     \be   (T-C)E(0) \leq 0.  \ee
The assertion (\ref{Appendix.4}) holds true. $\Box$

By a similar proof with Proposition (\ref{Appendix.1}), the following assertion holds.
 \begin{pro}\label{Appendix.5} There exists $T_1\ge0$ such that for any $T>T_1$, the only solution $u$ in $  C([0,T],H^1(\widehat{\Om}))$  to the system
\begin{equation}
\begin{cases}iu_t+\Delta_g u=0\qquad (x,t)\in \widehat{\Om}\times(0,T),
\cr  u=0\qquad(x,t)\in \omega\times
(0,T),\end{cases}
\end{equation}
is the trivial one $u\equiv 0$. \end{pro}


\end{document}